\documentclass[10pt]{article}
\usepackage{amssymb,amsthm,amsmath,amsthm}
\usepackage[numbers,sort&compress]{natbib}
\newtheorem{theorem.}{\textbf{Theorem}}[section]
\newtheorem{lemma.}{\textbf{Lemma}}[section]
\newtheorem{remark.}{\textbf{Remark}}
\newtheorem{definition.}{\textbf{Definition}}
\makeatletter
\def\blfootnote{\xdef\@thefnmark{}\@footnotetext}
\makeatother
\date{}
\begin{document}

\title{\Large{\bf Solutions for fourth-order
Kirchhoff type elliptic equations involving concave-convex
nonlinearities in $\mathbb{R}^{N}$}
\thanks{This work is supported by NSF of China (No.11801472) and NSF of China (No.11701463) and China Scholarship
Council (No.201708515186) and Science and Technology Innovation Team
of Education Department of Sichuan for Dynamical System and its
Applications (No.18TD0013), Youth Science and Technology Innovation
Team of Southwest Petroleum University for Nonlinear Systems
(No.2017CXTD02).}}

\author{Dong-Lun Wu$^{a,b}$\footnote{Corresponding author.}, Fengying Li$^{c}$\\
{\small $^{a}$College of Science, Southwest Petroleum University,}\\
{\small Chengdu, Sichuan 610500, P.R. China}\\
{\small $^{b}$Institute of Nonlinear Dynamics, Southwest Petroleum University,}\\
{\small Chengdu, Sichuan 610500, P.R. China}\\
{\small $^{c}$The School of Economic and Mathematics, Southwestern University of Finance and Economics, }\\
{\small Chengdu 611130, P. R. China}}
\date{}
\maketitle

\blfootnote{E-mail: wudl2008@163.com}

{\bf Abstract} In this paper, we show the existence and multiplicity
of solutions for the following fourth-order Kirchhoff type elliptic
equations
\begin{eqnarray*}
\Delta^{2}u-M(\|\nabla u\|_{2}^{2})\Delta u+V(x)u=f(x,u),\ \ \ \ \ x\in \mathbb{R}^{N},
\end{eqnarray*}
where $M(t):\mathbb{R}\rightarrow\mathbb{R}$ is the Kirchhoff
function, $f(x,u)=\lambda k(x,u)+ h(x,u)$, $\lambda\geq0$, $k(x,u)$ is
of sublinear growth and $h(x,u)$ satisfies some general 3-superlinear
growth conditions at infinity. We show the existence of at least one
solution for above equations for $\lambda=0$. For $\lambda>0$ small
enough, we obtain at least two nontrivial solutions. Furthermore, if
$f(x,u)$ is odd in $u$, we show that above equations possess
infinitely many solutions for all $\lambda\geq0$. Our theorems
generalize some known results in the literatures even for $\lambda=0$
and our proof is based on the variational methods.


\section{Introduction}

In this paper, we study the existence of multiple solutions for the
following fourth-order Kirchhoff type elliptic equations
\begin{eqnarray}
\Delta^{2}u-M(\|\nabla u\|_{2}^{2})\Delta u+V(x)u=f(x,u),\ \ \ \ \ x\in \mathbb{R}^{N},\label{1}
\end{eqnarray}
where $M\in C([0,+\infty),\mathbb{R})$ is a Kirchhoff-type function,
potential $V(x)\in C(\mathbb{R},\mathbb{R}^{+})$ and nonlinearities
$f(x,u)\in C(\mathbb{R}^{N}\times \mathbb{R},\mathbb{R})$, $1<N<8$.
The Kirchhoff type problems  on a bounded domain is introduced as
\begin{eqnarray}
\left\{
\begin{array}{ll}
-\left(a+b\displaystyle\int_{\Omega}|\nabla u|^{2}dx\right)\Delta
u=f(x,u)&\mbox{in $\Omega$},\\
u(x)=0,&\mbox{on $\partial\Omega$},
\end{array}
\right.\label{4}
\end{eqnarray}
which is related to the stationary analogue of the Kirchhoff equation
\begin{eqnarray}
u_{tt}-\left(a+b\int_{\Omega}|\nabla u|^{2}dx\right)\Delta
u=g(x,u).\label{2}
\end{eqnarray}
Equation (\ref{2}) was proposed by Kirchhoff in 1883 as a
generalization of following d'Alembert's wave equation
\begin{eqnarray*}
\rho\frac{\partial^{2}u}{\partial^{2}t}-\left(\frac{P_{0}}{h}+\frac{E}{2L}\int_{0}^{L}|\frac{\partial
u}{\partial x}|^{2}dx\right)\frac{\partial^{2} u}{\partial^{2}
x}=g(x,u)
\end{eqnarray*}
for free vibrations of elastic strings. It well known that, as a
useful model, the Kirchhoff equation has many applications in
mechanical and biological problems. After the work of Lions \cite{32},
the existence and multiplicity of solutions for Kirchhoff equations
have been studied by many mathematicians via the variational methods.
In recent years, some authors considered the Kirchhoff equations or
the $p$-Kirchhoff equations with concave-convex nonlinearities. We
remained the readers with references
\cite{36,40,37,33,39,LL,LLZ,23,38,35}. However, the nonlinearities
were required to satisfy some specific form in these papers, such as
$$f(x,u)=\varpi s_{1}(x)|u|^{\tau_{1}-2}u+s_{2}(x)|u|^{\tau_{2}-2}u.$$
In the present paper, we study the concave-convex nonlinearities with
abstract forms.

Subsequently, we recall some known results about the fourth-order
Kirchhoff type elliptic equations which has been studied by many
mathematicians \cite{1,20,4,6,7,21,10,12,13}. In 2012, Wang and An
\cite{9} considered problem (\ref{1}) in a bounded domain with
potential $V(x)\equiv0$ and $M(t)\equiv const.$ when $t$ large enough.
Under some 2-superlinear conditions, the authors obtain a nonnegative
solution by using the Mountain Pass Theorem. This problem is related
to the stationary analog of the evolution equation of Kirchhoff type
$$ u_{tt} + \Delta^2 u - \left( a + b \int_{\Omega} |\nabla u|^2 dx \right) \Delta u = f(x,u). $$
Actually, if $M(t)=a+bt$, we can obtain solutions for problem
(\ref{1}) when the growth of the nonlinearities is required to be
4-superlinear which has been shown by some previous works. Whether
there are solutions for problem (\ref{1}) with 2-superlinear
nonlinearities when $M(t)=a+bt$ is still open. In this paper, we only
consider the 4-superlinear case. In order to study different
4-superlinear nonlinearities, some different kinds of growth
conditions were introduced. In \cite{21}, Song and Chen showed problem
(\ref{1}) possesses infinitely many solutions under the following
monotonous condition.

$(MC)$ There exists $\nu\geq1$ such that
\begin{eqnarray*}
\nu\widetilde{F}(x,t)\geq\widetilde{F}(x,s t)
\end{eqnarray*}
for all $(x,t)\in \mathbb{R}^{N}\times \mathbb{R}$ and $s\in[0,1]$,
where $\widetilde{F}(x,t)=\frac{1}{4}f(x,t)t-F(x,t)$ and $F(x,t) =
    \int_0^t f(x,s)ds$.

By replacing $(MC)$ with  the following local $(AR)$-type condition,
Song and Chen \cite{21} also obtained infinitely many solutions for
problem (\ref{1}).

$(AR)$ There exist $l_{1} > 0$ and $\mu>4$ such that
$$f(x,t)t\geq\mu F(x,t) , \quad \text{ for a.e.~}x \in \mathbb{R}^N \text{ and }\forall\ |t| \geq l_{1}. $$

In 2015, Xu and Chen \cite{13} considered problem (\ref{1}) in
$\mathbb{R}^{3}$ and obtained infinitely many solutions under the
following superlinear condition which is weaker than $(AR)$.

$(WAR)$ There exist constants $s_{0}$, $l_{2}>0$ and
$\iota>\frac{3}{2}$ such that
\begin{eqnarray*}
s_{0}|t|^{2\iota}\widetilde{F}(x,t)\geq |F(x,t)|^{\iota}\quad \text{ for a.e.~}x \in \mathbb{R}^N \text{ and }\forall\ |t| \geq l_{2}.
\end{eqnarray*}
Obviously, we can obtain the following condition with $(MC)$, $(AR)$
and $(WAR)$ respectively.

$(FSL)$ $\widetilde{F}(x,t)\geq0$ for a.e. $x \in \mathbb{R}^N $ and $
|t| $ large enough.

In a recent paper, Ding and Li \cite{20} considered a class of
nonhomogenous fourth-order Kirchhoff equations with $(FSL)$. They
obtained the following theorem.

\begin{theorem.}\label{thDL}(See \cite{20})
Assume that $f(x,t)=w(x,t)+g(x)$, where $g \in L^2 (\mathbb{R}^N)$, $g
\not\equiv 0$ and the following conditions hold.

\begin{itemize}
\item[$(V')$] $V \in C(\mathbb{R}^N , \mathbb{R})$ satisfies
    $\displaystyle\inf_{x \in \mathbb{R}^N} V (x) \geq V_0 > 0$ and
    for each $M
    > 0,$  meas$\{x \in \mathbb{R}^N : V (x) \leq M \}<+\infty$,
    where $V_0$ is a constant and meas denotes the Lebesgue measure
    in $\mathbb{R}^N$.
\item[$(w_{1})$] $w \in C(\mathbb{R}^N\times\mathbb{R},\mathbb{R} )$
    and
$$
|w(x,t)| \leq C(1+|t|^{p-1}) \quad \text{for some } 4 < p < 2^*=
\begin{cases}
\frac{2N}{N-4}, & 8 > N >4\\
+ \infty, & 1 < N \leq 4,
\end{cases}
$$
where $C$ is a positive constant.
\item[$(w_{2})$] $w(x,t) = o(|t|)$ as $|t| \to 0$ uniformly in $x
    \in \mathbb{R}^N$.
\item[$(w_{3})$] $\frac{W(x,t)}{t^4} \to + \infty$ as $|t| \to +
    \infty$ uniformly in $x \in \mathbb{R}^N$, \text{where $W(x,t) =
    \int_0^t w(x,s)ds$.}.
\item[$(w_{4})$] There exist $L > 0$ and $d \in \left[ 0,
    \frac{V_0}{2} \right]$ such that
$$ \frac{1}{4}w(x,t)t-W(x,t) \geq -\frac{d}{4}|t|^2, \quad \text{ for a.e.~}x \in \mathbb{R}^N \text{ and }\forall |t| \geq L. $$
\end{itemize}

Then, there exists a constant $g_0 > 0$ such that the problem
(\ref{1}) has at least two different solutions whenever $\|g\|_{L^2} <
g_0$ and $M(t)=a+bt$ with $a,b>0$, one is negative energy solution,
and the other is positive energy solution.
\end{theorem.}

\begin{remark.}
The condition $(w_{4})$ is required to hold for any $t\geq L$ in
\cite{20}. However, the authors used this condition for any $|t|\geq
L$ implicitly.
\end{remark.}

In order to use the variational methods to obtain the results, it is
not enough to show the geometric structure of the corresponding
functional. We also need to guarantee the convergence of the
asymptotic critical sequence which can be obtained by the compactness
of the embedding. However, since the domain is unboundedness, there is
no natural compact embedding to use. To overcome this difficulty,
periodic, coercive and radial symmetric conditions are put forward.
Condition $(V')$ is a classical coercive condition on $V(x)$ to make
sure the embedding is compact. It has been shown by Bartsch and Wang
in \cite{2} that the following coercive condition is weaker than
$(V')$.

$(V)$\ $V\in C(\mathbb{R}^{N},\mathbb{R})$,
$\inf_{x\in\mathbb{R}^{N}}V(x)>0$. There exists $\bar{r}>0$ such that
\begin{eqnarray*}
\lim_{|y|\rightarrow\infty}meas\left\{x\in \mathbb{R}^{N}: |x-y|\leq
\bar{r},V(x)\leq \mathcal{M}\right\}=0,\ \ \ \forall\ \mathcal{M}>0.
\end{eqnarray*}
There are still some other ways to get the compactness
back(see\cite{27,28,30,29}). In this paper, we consider the coercive
case and use condition $(V)$ to obtain the compactness of the
embedding.

Before we state our results, we introduce some conditions on $M$. In
problem (\ref{1}), the Kirchhoff function $M(t)$ is a abstract
function, which has been rarely considered when the problem lies in
$\mathbb{R}^{N}$. Through out this paper, we assume that $M(t)$
satisfies the following conditions.

$(M_{1})$ $M \in C([0,+\infty),\mathbb{R})$ and there exists $m_{0}>0$
such that $M(t) \geq m_{0}$ for all $t \in [0,+\infty)$.

$(M_{2})$ There exist positive constants $\sigma_{1}$ and $\sigma_{2}$
such that
$$\sigma_{2}(t^{2}+t)\geq\widehat{M}(t)\geq \frac{1}{2}
M(t)t+\sigma_{1} t\ \ \ \ \mbox{for any}\ \  t\geq0,$$ where
$\widehat{M}(t)=\int_{0}^{t}M(s)ds$.

\begin{remark.}
It is easy to check that the original Kirchhoff function $M(t)=a+bt$
for any $a,b>0$ satisfies $(M_{1})$-$(M_{2})$. There are still some
other functions admitting our conditions, such as
\begin{eqnarray*}\widehat{M}(t)=\frac{at+bt^{2}}{1+\ln(t^{2}+1)}\ \ \mbox{with $a, b>0$}.
\end{eqnarray*}
\end{remark.}
Subsequently, in order to study the concave-convex nonlinearities, we
consider
\begin{eqnarray}
f(x,t)=\lambda k(x,t)+ h(x,t).\label{5}
\end{eqnarray}
Letting $K(x,t)=\int^{t}_{0}k(x,v)dv$ and
$H(x,t)=\int^{t}_{0}h(x,v)dv$, we state our main theorems.

\begin{theorem.}\label{th1.3} Suppose that
(\ref{5}), $(V)$, $(M_{1})$, $(M_{2})$ and the following conditions
hold

\begin{itemize}
\item[$(f_{1})$] there exit $\bar{x}\in \mathbb{R}^{N}$,
    $r_{0}\in(1,2)$ and $b_{0}>0$ such that $K(\bar{x},t)\geq
    b_{0}|t|^{r_{0}}$ for all $t\in \mathbb{R}$;

\item[$(f_{2})$] for any $(x,t)\in \mathbb{R}^{N}\times \mathbb{R}$,
    there exist $r_{1}$, $r_{2}\in(1,2)$ such that
\begin{eqnarray*}
|k(x,t)|\leq b_{1}(x)|t|^{r_{1}-1}+b_{2}(x)|t|^{r_{2}-1},
\end{eqnarray*}
where $b_{i}(x)\in L^{\beta_{i}}(\mathbb{R}^{N},\mathbb{R}^{+})$ and
$\beta_{i}\in\left(\frac{2^{*}}{2^{*}-r_{i}},\frac{2}{2-r_{i}}\right]$
for $i=1,2$;

\item[$(f_{3})$]  $h(x,t)=o(|t|)$ as $t\rightarrow0$ uniformly in $x
    \in \mathbb{R}^N$;

\item[$(f_{4})$] $H(x,t)/t^{4}\rightarrow+\infty$ as
    $|t|\rightarrow\infty$ uniformly in $x \in \mathbb{R}^N$;

\item[$(f_{5})$] there exist positive constants $d_{1}$ and
    $\rho_{\infty}$ such that
\begin{eqnarray*}
\widetilde{H}(x,t)=\frac{1}{4}h(x,t)t- H(x,t)\geq -d_{1}t^{2}\ \ \ \mbox{for all}\ \ \text{ a.e.~}x \in \mathbb{R}^N\ \ \mbox{and}\ \ \ |t|\geq\rho_{\infty};
\end{eqnarray*}

\item[$(f_{6})$] there exist $d_{2}>0$ and $4<\zeta<2^{*}$ such that
\begin{eqnarray*}
|h(x,t)|\leq d_{2}(|t|+|t|^{\zeta-1})\ \ \ \mbox{for all}\ \ \text{ a.e.~}x \in \mathbb{R}^N\ \ \mbox{and}\ \
t\in  \mathbb{R}.
\end{eqnarray*}
\end{itemize}

Then problem (\ref{1}) possesses at least one solution for $\lambda=0$
and there exists $\lambda_{1}>0$ such that for any
$\lambda\in(0,\lambda_{1})$, problem (\ref{1}) possesses at least two
solutions.
\end{theorem.}

\begin{theorem.}\label{th1.4} Suppose that
(\ref{5}), $(V)$, $(M_{1})$, $(M_{2})$, $(f_{1})$-$(f_{6})$ and

\noindent$(f_{7})$ $K(x,-t)=K(x,t)$ and $H(x,-t)=H(x,t)$ for all
$(x,t)\in \mathbb{R}^{N}\times \mathbb{R}$.

Then for any $\lambda\geq0$, problem (\ref{1}) possesses infinitely
many solutions.
\end{theorem.}

\begin{remark.} Obviously, condition $(f_{5})$ is weaker than $(w_{4})$. We can also see
that $(f_{5})$ is weaker than $(MC)$, $(AR)$ and $(WAR)$. Hence
Theorem \ref{th1.3} and \ref{th1.4} generalize  Theorems 1.1, 1.2, 1.3
in \cite{21}  and Theorems 1.1, 1.2 in \cite{13}.
\end{remark.}
\begin{remark.} In
Theorems \ref{th1.3} and \ref{th1.4}, the sign of $H(x,t)$ is
indefinite. Although we have $(f_{4})$, $H(x,t)$ can also be negative
around origin with respect to $t$.
\end{remark.}
\begin{remark.}
Although there were some papers concerning on the fourth-order
Kirchhoff type elliptic equations with concave-convex nonlinearities
on bounded domain, to the best of the knowledge of the authors, this
is the first work on fourth-order Kirchhoff type elliptic equations
with concave-convex growth on unbounded domain.
\end{remark.}

In this paper, we will use the variational method to prove our
theorems. First, we introduce the definition of the $(PS)^{*}$
condition.

\begin{definition.} Let $E$ be a Hilbert space. A functional $I\in
C^{1}(E,R)$ is said to satisfy the $(PS)^{*}$ condition with respect
to $E_{j}$, $j=1,2,\cdots$, if any sequence $x_{j}\in E_{j}$
satisfying
\begin{eqnarray*}
|I(x_{j})|<\infty\ \ \ \mbox{and}\ \
I'|_{E_{j}}(x_{j})\rightarrow0,
\end{eqnarray*}
imply a convergent subsequence, where $E_{j}$ is a sequence of linear
subspace of $E$ with finite dimensional.
\end{definition.}

The following critical point theorem is needed to obtain the
multiplicity of solutions.

\begin{lemma.} (Chang\cite{22})\label{leZhang} Suppose that
 $E$ is a Hilbert space, $I\in C^{1}(E,\mathbb{R})$ is even with $I(0)=0$, and
that
\begin{itemize}
\item[$(C_{1})$] there exist $\varrho$, $\alpha>0$ and a finite
    dimensional linear subspace $X$ such that $I|_{X^{\bot}\bigcap
\partial B_{\varrho}}\geq\alpha$, where $B_{\varrho}=\{u\in E: \|u\|_{E}\leq
\varrho\}$;

\item[$(C_{2})$] there is a sequence of linear subspaces
    $\tilde{X}_{m}$, $dim\tilde{X}_{m}=m$, and there exists
    $r_{m}>0$ such that
\begin{eqnarray*}
I(u)\leq0\ \ \mbox{on}\ \ \tilde{X}_{m}\setminus B_{r_{m}},\
m=1,2,\cdots.
\end{eqnarray*}
\end{itemize}

If, further, $I$ satisfies the $(PS)^{*}$ condition
 with respect to
$\{\tilde{X}_{m}| m=1,2,\cdots\}$, then $I$ possesses infinitely many
distinct critical points corresponding to positive critical values.
\end{lemma.}

\section{Preliminaries}

In this paper, we let
\begin{eqnarray*}
H^{2}(\mathbb{R}^{N})=\left\{u\in L^{2}(\mathbb{R}^{N}):\nabla u\in L^{2}(\mathbb{R}^{N}),\Delta u\in L^{2}(\mathbb{R}^{N})\right\}
\end{eqnarray*}
with the norm
\begin{eqnarray*}
\|u\|_{H^{2}}^{2}=\int_{\mathbb{R}^{N}}(|\Delta u|^{2}+|\nabla u|^{2}+u^{2})dx
\end{eqnarray*}
Set
\begin{eqnarray*}
E=\left\{u\in H^{2}(\mathbb{R}^{N}):\int_{\mathbb{R}^{N}}(|\Delta u|^{2}+|\nabla u|^{2}+V(x)u^{2})dx<+\infty\right\}
\end{eqnarray*}
with the inner product
\begin{eqnarray*}
\langle u,v\rangle_{E}=\int_{\mathbb{R}^{N}}(\Delta u\cdot\Delta v+\nabla u\cdot\nabla
v+V(x)uv)dx
\end{eqnarray*}
and the norm $\|u\|_{E}=\langle u,u\rangle^{1/2}$. Obviously, It is
well known that under hypothesis $(V)$, the embedding $E
\hookrightarrow L^{s} (\mathbb{R}^{N})$ is continuous for $s \in [2,
2^{*}]$ and compact for $s \in [2, 2^{*})$. Then, for any $s \in [2,
2^{*}]$, there exists $C_{s}>0$ such that
\begin{eqnarray}
\|u\|_{s}\leq C_{s}\|u\|_{E}\ \ \ \mbox{for all}\ \ u\in
E.\label{3}
\end{eqnarray}
It is known that the weak solutions for problem (\ref{1}) are the
critical points of the following functional
\begin{eqnarray*}
I(u)&=&\frac{1}{2}\|\Delta u\|_{2}^{2}+\frac{1}{2}\widehat{M}(\|\nabla u\|_{2}^{2})+\frac{1}{2}\int_{\mathbb{R}^{N}}V(x)u^{2}dx-\int_{\mathbb{R}^{N}}F(x,u)dx\\
&=&\frac{1}{2}\|\Delta u\|_{2}^{2}+\frac{1}{2}\widehat{M}(\|\nabla u\|_{2}^{2})+\frac{1}{2}\int_{\mathbb{R}^{N}}V(x)u^{2}dx-\lambda\int_{\mathbb{R}^{N}}K(x,u)dx-\int_{\mathbb{R}^{N}}H(x,u)dx.
\end{eqnarray*}
Similar to the proof of Proposition 2.2 in \cite{29}, under
 $(f_{2})$ and $(f_{3})$, we see that $I\in
C^{1}(E,\mathbb{R})$ and for each $u$, $v\in E$,
\begin{eqnarray}
\langle I'(u),v \rangle&=&\int_{\mathbb{R}^{N}}\Delta u\cdot \Delta vdx+M(\|\nabla u\|_{2}^{2})\int_{\mathbb{R}^{N}} \nabla u\cdot\nabla
vdx\nonumber\\
&&+\int_{\mathbb{R}^{N}}V(x)uvdx-\lambda\int_{\mathbb{R}^{N}}k(x,u)vdx-\int_{\mathbb{R}^{N}}h(x,u)vdx.\label{11}
\end{eqnarray}

\section{Proof of Theorem \ref{th1.3}}

\begin{lemma.}\label{le3.1} Assume (\ref{5}), $(V)$, $(M_{1})$, $(M_{2})$, $(f_{2})$, $(f_{3})$ and $(f_{6})$ hold,
then there exists $\lambda_{1}>0$ such that for all
$\lambda\in[0,\lambda_{1})$, there exist $\varrho$, $\alpha>0$ such
that $I|_{\partial B_{\varrho}}\geq\alpha$ , where $B_{\varrho}=\{u\in
E: \|u\|_{E}\leq \varrho\}$.
\end{lemma.}

{\bf Proof.} By $(f_{2})$, we obtain that
\begin{eqnarray}
|K(t,x)|\leq
\frac{1}{r_{1}}b_{1}(x)|t|^{r_{1}}+\frac{1}{r_{2}}b_{2}(x)|t|^{r_{2}}\label{27}
\end{eqnarray}
for all $(x,t)\in \mathbb{R}^{N}\times \mathbb{R}$. It follows from
$(f_{3})$ and $(f_{6})$, for any $\varepsilon>0$, there exists
$D_{\varepsilon}>0$ such that
\begin{eqnarray}
|H(x,t)|\leq\varepsilon|t|^{2}+D_{\varepsilon}|t|^{\zeta},\ \ \ \ \
 \forall (x,t) \in \mathbb{R}^{N}\times\mathbb{R}.\label{42}
\end{eqnarray}
By (\ref{3}), $(M_{1})$, $(M_{2})$, (\ref{27}) and (\ref{42}), we have
\begin{eqnarray*}
&& I(u)\\
&=&
\frac{1}{2}\|\Delta u\|_{2}^{2}+\frac{1}{2}\widehat{M}(\|\nabla u\|_{2}^{2})+\frac{1}{2}\int_{\mathbb{R}^{N}}V(x)u^{2}dx-\lambda\int_{\mathbb{R}^{N}}K(x,u)dx-\int_{\mathbb{R}^{N}}H(x,u)dx\nonumber\\
&\geq& \mbox{}
\frac{1}{2}\|\Delta u\|_{2}^{2}+\left(\frac{m_{0}}{2}+\sigma_{1}\right)\|\nabla u\|_{2}^{2}+\frac{1}{2}\int_{\mathbb{R}^{N}}V(x)u^{2}dx\\
&&-\lambda\left(\frac{1}{r_{1}}\int_{\mathbb{R}^{N}}b_{1}(x)|u|^{r_{1}}dx+\frac{1}{r_{2}}\int_{\mathbb{R}^{N}}b_{2}(x)|u|^{r_{2}}dx
\right)-\left(\varepsilon\int_{\mathbb{R}^{N}}|u|^{2}dx+D_{\varepsilon}\int_{\mathbb{R}^{N}}|u|^{\zeta}dx\right)\nonumber\\
&\geq& \mbox{}
\min\left\{\frac{1}{2},\frac{m_{0}}{2}+\sigma_{1}\right\}\|u\|_{E}^{2}-\lambda M_{1}(
\|u\|_{E}^{r_{1}-2}+\|u\|_{E}^{r_{2}-2})-\left(\varepsilon C_{2}^{2}\|u\|_{E}^{2}+D_{\varepsilon}C_{\zeta}^{\zeta}\|u\|_{E}^{\zeta}\right)\nonumber\\
&\geq&\mbox{} \left(\min\left\{\frac{1}{2},\frac{m_{0}}{2}+\sigma_{1}\right\}-\lambda M_{1}(
\|u\|_{E}^{r_{1}-2}+\|u\|_{E}^{r_{2}-2})-\varepsilon C_{2}^{2}-D_{\varepsilon}C_{\zeta}^{\zeta}\|u\|_{E}^{\zeta-2}\right)\|u\|_{E}^{2},\label{4}
\end{eqnarray*}
where
$M_{1}=\max\left\{\frac{1}{r_{1}}C_{r_{1}\beta_{1}^{*}}^{r_{1}}\|b_{1}\|_{\beta_{1}},\frac{1}{r_{2}}C_{r_{2}\beta^{*}}^{r_{2}}\|b_{2}\|_{\beta}\right\}
$. If we choose $\varepsilon$ small enough, it is easy to see that
there exist positive constants $\lambda_{1}$, $\varrho$ and $\alpha$
such that $I|_{\partial B_{\varrho}}\geq\alpha$ for all
$\lambda\in[0,\lambda_{1})$. We finish the proof of this lemma.

\begin{lemma.}
Suppose (\ref{5}), $(V)$, $(M_{2})$, $(f_{2})$ and $(f_{4})$
 hold, then there exists $\tilde{e} \in E$ such that
$\|\tilde{e}\|> \varrho$ and $I(\tilde{e}) \leq 0$, where $\varrho$ is
defined in Lemma \ref{le3.1}. \end{lemma.}

{\bf Proof.} Choose $e \in C^{\infty}_{0}(\Upsilon_{1}(0),\mathbb{R})$
such that $\|e\|_{E}=1$, where
$\Upsilon_{r}(x_{0})=\{x\in\mathbb{R}^{N}: |x-x_{0}|\leq r\}$.  We can
see that there exist $L_{0}>0$ and $\Sigma\subset \Upsilon_{1}(0)$
such that $|e|\geq L_{0}$ for all $x\in\Sigma$ with $meas (\Sigma)>0$.
By $(f_{4})$, for any $A>0$ there exists $Q>0$ such that
\begin{eqnarray*}
\frac{H(x,t)}{t^{4}}\geq A
\end{eqnarray*}
for all $|t|\geq Q$ and $x\in\mathbb{R}^{N}$, which implies that
\begin{eqnarray}
\int_{\Sigma}\frac{H(x,\xi e)}{|\xi e|^{4}}dx\geq A
meas (\Sigma),\label{10}
\end{eqnarray}
for all $\xi\geq Q/L_{0}$. By~(\ref{3}), $(M_{2})$, (\ref{27}) and
(\ref{10}), for any $\xi>0$ large enough, we have
\begin{eqnarray*}
\frac{I(\xi e)}{\xi^{4}}
&=&\frac{1}{2\xi^{2}}\|\Delta e\|_{2}^{2}+\frac{1}{2\xi^{4}}\widehat{M}(\xi^{2}\|\nabla e\|_{2}^{2})+\frac{1}{2\xi^{2}}\int_{\Sigma}V(x)e^{2}dx\\
&&-\frac{\lambda}{\xi^{4}}\int_{\Sigma}K(x,\xi e)dx-\frac{1}{\xi^{4}}\int_{\Sigma}H(x,\xi e)dx\nonumber\\
&\leq& \mbox{}
\frac{1}{2\xi^{2}}\|\Delta e\|_{2}^{2}+\frac{\sigma_{2}}{2\xi^{4}}\left(\xi^{4}\|\nabla e\|_{2}^{4}+\xi^{2}\|\nabla e\|_{2}^{2}\right)+\frac{1}{2\xi^{2}}\int_{\Sigma}V(x)e^{2}dx\\
&&+\lambda
M_{1}(\xi^{r_{1}-4}\|e\|_{E}^{r_{1}}+\xi^{r_{2}-4}\|e\|_{E}^{r_{2}})-L_{0}^{4}\int_{\Sigma}\frac{H(x,\xi
e)}{|\xi e|^{4}}dx\nonumber\\
&\leq& 1+\frac{\sigma_{2}}{2}-L_{0}^{4}Ameas (\Sigma).
\end{eqnarray*}
By the arbitrariness of $A$, there exists $\xi_{0}>0$ such that
$I(\xi_{0}e)<0$ and $\|\xi_{0}e\|> \varrho$. Let $\tilde{e}=\xi_{0}e$,
we can see $I(\tilde{e})<0$, which proves this lemma.

\begin{lemma.}\label{le3.3}
Suppose (\ref{5}), $(V)$, $(M_{1})$, $(M_{2})$ and $(f_{2})$-$(f_{6})$
 hold, then $I$ satisfies the $(PS)$ condition.
\end{lemma.}

{\bf Proof.} Let $\{u_{n}\} \subset E$ be a sequence such that
$\{I(u_{n})\}$ is bounded and $I^{\prime}(u_{n}) \rightarrow 0$ as $n
\rightarrow \infty$. Then there exists a constant $M_{2} > 0$ such
that
\begin{eqnarray}
\ |I(u_{n})| \leq M_{2},\ \ \ \ \ \|I^{\prime}(u_{n})\|_{E^{*}}
\leq M_{2}.\label{32}
\end{eqnarray}
Next, we show that $\{u_{n}\}$ is bounded in $E$. Arguing in an
indirect way, we assume that $\|u_{n}\|_{E}\rightarrow+\infty$ as
$n\rightarrow \infty$. Set $z_{n}=\frac{u_{n}}{\|u_{n}\|_{E}}$, then
$\|z_{n}\|_{E}=1$, which implies that there exists a subsequence of
$\{z_{n}\}$, still denoted by $\{z_{n}\}$, such that
$z_{n}\rightharpoonup z_{0}$ in $E$ and $z_{n}\rightarrow z_{0}$ a.e.
in $\mathbb{R}^{N}$ as $n\rightarrow \infty$. From $(f_{3})$, we can
deduce that $\widetilde{H}(x,t)=o(t^{2})$ as $|t|\rightarrow0$
uniformly in $x$, then there exists $\rho_{0}\in(0,\rho_{\infty})$
such that
\begin{eqnarray}
|\widetilde{H}(x,t)|\leq t^{2}\label{14}
\end{eqnarray}
for all $|t|\leq\rho_{0}$ and $x\in\mathbb{R}^{N}$ . If
$z_{0}\equiv0$, we can deduce
\begin{eqnarray*}
&&o(1){\nonumber}\\&=& \frac{
M_{2}+\frac{1}{4}M_{2}\|u_{n}\|_{E}}{\|u_{n}\|_{E}^{2}}\\ &\geq&
\frac{ I(u_{n})-\frac{1}{4}\langle I^{\prime}(u_{n}),u_{n}\rangle}{\|u_{n}\|_{E}^{2}}{\nonumber}\\
&\geq& \mbox{}
\frac{1}{\|u_{n}\|_{E}^{2}}\left(\frac{1}{4}\|\Delta u_{n}\|_{2}^{2}+\frac{1}{2}\widehat{M}(\|\nabla u_{n}\|_{2}^{2})-\frac{1}{4}M(\|\nabla u_{n}\|_{2}^{2})\|\nabla u_{n}\|_{2}^{2}+\frac{1}{4}\int_{\mathbb{R}^{N}}V(x)u_{n}^{2}dx\right){\nonumber}\\
& & \mbox{}-\frac{2\lambda M_{1}}{\|u_{n}\|_{E}^{2}}\left(\|u_{n}\|_{E}^{r_{1}}+\|u_{n}\|_{E}^{r_{2}}\right)+\frac{1}{\|u_{n}\|_{E}^{2}}\int_{\mathbb{R}^{N}}\widetilde{H}(x,u_{n})dx{\nonumber}\\
&\geq&\frac{1}{\|u_{n}\|_{E}^{2}}\left(\frac{1}{4}\|\Delta u_{n}\|_{2}^{2}+\frac{\sigma_{1}}{2}\|\nabla u_{n}\|_{2}^{2}+\frac{1}{4}\int_{\mathbb{R}^{N}}V(x)u_{n}^{2}dx\right)-\frac{2\lambda M_{1}}{\|u_{n}\|_{E}^{2}}\left(\|u_{n}\|_{E}^{r_{1}}+\|u_{n}\|_{E}^{r_{2}}\right){\nonumber}\\
& & +\frac{1}{\|u_{n}\|_{E}^{2}}\left(\int_{|u_{n}|\leq \rho_{0}}\widetilde{H}(x,u_{n})dx+\int_{|u_{n}|>\rho_{\infty}}\widetilde{H}(x,u_{n})dx\right)+\frac{1}{\|u_{n}\|_{E}^{2}}\int_{\rho_{0}<|u_{n}|\leq\rho_{\infty}}\widetilde{H}(x,u_{n})dx{\nonumber}\\
&\geq& \mbox{}
\min\left\{\frac{1}{4},\frac{\sigma_{1}}{2}\right\}-\frac{1}{\|u_{n}\|_{E}^{2}}\left(\int_{|u_{n}|\leq \rho_{0}}|u_{n}|^{2}dx+d_{1}\int_{|u_{n}|>\rho_{\infty}}|u_{n}|^{2}dx\right){\nonumber}\\
& & \mbox{}-\frac{(d_{2}+1+D_{1})(\rho_{\infty}^{2}+\rho_{\infty}^{\zeta})}{\rho_{0}^{2}}\int_{\rho_{0}<|u_{n}|\leq\rho_{\infty}}|z_{n}|^{2}dx+o(1){\nonumber}\\
&\geq& \mbox{}
\min\left\{\frac{1}{4},\frac{\sigma_{1}}{2}\right\}-\left(1+d_{1}+\frac{(d_{2}+1+D_{1})(\rho_{\infty}^{2}+\rho_{\infty}^{\zeta})}{\rho_{0}^{2}}\right)\int_{\mathbb{R}^{N}}|z_{n}|^{2}dx+o(1){\nonumber}\\
&\rightarrow&\min\left\{\frac{1}{4},\frac{\sigma_{1}}{2}\right\}\ \ \ \mbox{as}\ \
n\rightarrow\infty,
\end{eqnarray*}
which is a contradiction. Then we have $z_{0}\not\equiv0$. Let
$\Omega=\{x\in \mathbb{R}^{N}|\ |z_{0}(x)|>0 \}$. Then we can see that
$meas(\Omega)>0$.  Since $\|u_{n}\|_{E}\rightarrow+\infty$ as
$n\rightarrow \infty$ and $|u_{n}|=|z_{n}|\cdot\|u_{n}\|_{E}$, then we
have $|u_{n}|\rightarrow+\infty$ as $n\rightarrow \infty$ for a.e.
$x\in\Omega$. It follows from $(f_{3})$ and $(f_{4})$ that there
exists
 $M_{3}>0$ such that
\begin{eqnarray}
H(x,t)\geq -M_{3}t^{2}\label{6}
\end{eqnarray}
for all $(x,t)\in \mathbb{R}\times\mathbb{R}^{N}$. Hence, by
(\ref{3}), we obtain
\begin{eqnarray*}
\int_{\mathbb{R}^{N}\setminus\Omega}\frac{H(x,u_{n})}{\|u_{n}\|_{E}^{4}}dx\geq -M_{3}\int_{\mathbb{R}^{N}\setminus\Omega}\frac{u_{n}^{2}}{\|u_{n}\|_{E}^{4}}dx\geq -M_{3}\frac{\|u_{n}\|_{2}^{2}}{\|u_{n}\|_{E}^{4}}\geq -M_{3}\frac{C_{2}^{2}}{\|u_{n}\|_{E}^{2}},
\end{eqnarray*}
which implies
\begin{eqnarray}
\liminf_{n\rightarrow\infty}\int_{\mathbb{R}^{N}\setminus\Omega}\frac{H(x,u_{n})}{\|u_{n}\|_{E}^{4}}dx\geq 0.
\end{eqnarray}
Moreover, we deduce from $(f_{4})$ and Fatou's Lemma that
\begin{eqnarray}
\liminf_{n\rightarrow
\infty}\int_{\Omega}\frac{H(x,u_{n})}{|u_{n}|^{4}}|z_{n}|^{4}dx\rightarrow+\infty\ \ \ \mbox{as}\ \ n\rightarrow\infty.
\end{eqnarray}
It follows from $(M_{2})$ that
\begin{eqnarray*}
&&2I(u_{n})+2\lambda\int_{\mathbb{R}^{N}}K(x,u_{n})dx+2\int_{\mathbb{R}^{N}}H(x,u_{n})dx\\
&=&
\|\Delta u_{n}\|_{2}^{2}+\widehat{M}(\|\nabla u_{n}\|_{2}^{2})+\int_{\mathbb{R}^{N}}V(x)u_{n}^{2}dx\nonumber\\
&\leq&\|\Delta u_{n}\|_{2}^{2}+\sigma_{2}(\|\nabla u_{n}\|_{2}^{4}+\|\nabla u_{n}\|_{2}^{2})+\int_{\mathbb{R}^{N}}V(x)u_{n}^{2}dx\nonumber\\
&\leq&(\sigma_{2}+1)\|u_{n}\|_{E}^{2}+\sigma_{2}\| u_{n}\|_{E}^{4},
\end{eqnarray*}
which implies that
\begin{eqnarray*}
\sigma_{2}&=&\liminf_{n\rightarrow\infty}\frac{(\sigma_{2}+1)\|u_{n}\|_{E}^{2}+\sigma_{2}\| u_{n}\|_{E}^{4}}{\| u_{n}\|_{E}^{4}}\\
&\geq&\liminf_{n\rightarrow\infty}\frac{2}{\| u_{n}\|_{E}^{4}}\left(I(u_{n})+\lambda\int_{\mathbb{R}^{N}}K(x,u_{n})dx+\int_{\mathbb{R}^{N}}H(x,u_{n})dx\right)\\
&\geq&\liminf_{n\rightarrow\infty}\frac{2}{\| u_{n}\|_{E}^{4}}\left(-\lambda
M_{1}\left(\|u_{n}\|_{E}^{r_{1}}+\|u_{n}\|_{E}^{r_{2}}\right)+\int_{\Omega}H(x,u_{n})dx\right)\\
&\geq&\liminf_{n\rightarrow\infty}2\int_{\Omega}\frac{H(x,u_{n})}{| u_{n}|^{4}}| z_{n}|^{4}dx\\
&\rightarrow&+\infty\ \ \ \mbox{as}\ \
n\rightarrow\infty,
\end{eqnarray*}
which is a contradiction. Hence $\{u_{n}\}$ is bounded in $E$. Then
there exists a subsequence, still denoted by $\{u_{n}\}$, such that
$u_{n}\rightharpoonup u$ in $E$. Therefore
\begin{eqnarray*}
\langle I^{\prime}(u_{n})-I^{\prime}(u),u_{n}-u\rangle\rightarrow0\
\ \ \mbox{as}\ \ \ n\rightarrow+\infty.
\end{eqnarray*}
By (\ref{3}) and $(f_{2})$, we have
\begin{eqnarray*}
&&\left|\int_{\mathbb{R}^{N}}(h(x,u_{n})-h(x,u),u_{n}-u)dx\right|\\
&\leq&\int_{\mathbb{R}^{N}}(|h(x,u_{n})|+|h(x,u)|)|u_{n}-u|dx{\nonumber}\\
&\leq& d_{2}\int_{\mathbb{R}^{N}}\left(|u_{n}|+|u_{n}|^{\zeta-1}+|u|+|u|^{\zeta-1}\right)|u_{n}-u|dx{\nonumber}\\
&\leq&
d_{2}\left(\|u_{n}\|_{2}\|u_{n}-u\|_{2}+\|u_{n}\|_{\zeta}^{\zeta-1}\|u_{n}-u\|_{\zeta}+\|u\|_{2}\|u_{n}-u\|_{2}+\|u\|_{\zeta}^{\zeta-1}\|u_{n}-u\|_{\zeta}\right)\\
&\leq&
d_{2}\left(C_{2}+C_{\zeta}^{\zeta-1}\right)\left(\|u_{n}\|_{E}\|u_{n}-u\|_{2}+\|u_{n}\|_{E}^{\zeta-1}\|u_{n}-u\|_{\zeta}+\|u\|_{E}\|u_{n}-u\|_{2}+\|u\|_{E}^{\zeta-1}\|u_{n}-u\|_{\zeta}\right)\\
&\rightarrow&0\ \ \ \mbox{as}\ \ n\rightarrow\infty.
\end{eqnarray*}
For $i=1,2$, set
\begin{alignat*}{2}
&S_{i,1}=\left(\frac{2^{*}}{2^{*}-r_{i}},\frac{22^{*}}{2^{*}-2(r_{i}-1)}\right],\ \ \ \
&S_{i,2}=\left(\frac{22^{*}}{2^{*}(2-r_{i})+2^{*}-2},\frac{2}{2-r_{i}}\right].
\end{alignat*}
It is easy to see that
$\left(\frac{2^{*}}{2^{*}-r_{i}},\frac{2}{2-r_{i}}\right]=
S_{i,1}\bigcup S_{i,2}$ and $S_{i,1}\bigcap S_{i,2}\neq\emptyset$.
Moreover, let
\begin{eqnarray*} \xi_{i}=\left\{
\begin{array}{ll}
2^{*}&\mbox{if $\beta_{i}\in S_{i,1}$},\\
2&\mbox{if $\beta_{i}\in S_{i,2}\setminus S_{i,1}$}.
\end{array}
\right.
\end{eqnarray*}
By an easy computation, we deduce that there exists
$\eta_{i}\in[2,2^{*})$ such that
$\frac{1}{\beta_{i}}+\frac{r_{i}-1}{\xi_{i}}+\frac{1}{\eta_{i}}=1$. It
follows from (\ref{3}) and $(f_{2})$ that
\begin{eqnarray*}
&&
\int_{\mathbb{R}^{N}}(k(x,u_{n})-k(x,u),u_{n}-u)dx{\nonumber}\\
&\leq&\int_{\mathbb{R}^{N}}|k(x,u_{n})-k(x,u)||u_{n}-u|dx{\nonumber}\\
&\leq& \mbox{}\int_{\mathbb{R}^{N}}(b_{1}(x)(|u_{n}|^{r_{1}-1}+|u|^{r_{1}-1})+b_{2}(x)(|u_{n}|^{r_{2}-1}+|u|^{r_{2}-1}))|u_{n}-u|dx{\nonumber}\\
&\leq&
\mbox{}\|b_{1}\|_{\beta_{1}}(\|u_{n}\|_{\xi_{1}}^{r_{1}-1}+\|u\|_{\xi_{1}}^{r_{1}-1})\|u_{n}-u\|_{\eta_{1}}+\|b_{2}\|_{\beta_{2}}(\|u_{n}\|_{\xi_{2}}^{r_{2}-1}+\|u\|_{\xi_{2}}^{r_{2}-1})\|u_{n}-u\|_{\eta_{2}}{\nonumber}\\
&\leq&
\mbox{}C_{\xi_{1}}^{r_{1}-1}\|b_{1}\|_{\beta_{1}}(\|u_{n}\|_{E}^{r_{1}-1}+\|u\|_{E}^{r_{1}-1})\|u_{n}-u\|_{\eta_{1}}+C_{\xi_{2}}^{r_{2}-1}\|b_{2}\|_{\beta_{2}}(\|u_{n}\|_{E}^{r_{2}-1}+\|u\|_{E}^{r_{2}-1})\|u_{n}-u\|_{\eta_{2}}{\nonumber}\\
&\rightarrow&0\ \ \ \mbox{as}\ \ n\rightarrow\infty.
\end{eqnarray*}
Therefore,
$\int_{\mathbb{R}^{N}}(f(x,u_{n})-f(x,u))(u_{n}-u)dx\rightarrow0$ as
$n\rightarrow\infty$. It follows from (\ref{11}) that
\begin{eqnarray*}
&&\langle
I^{\prime}(u_{n})-I^{\prime}(u),u_{n}-u\rangle\nonumber\\
&=&\|\Delta (u_{n}-u)\|_{2}^{2}+M(\|\nabla
u_{n}\|_{2}^{2})\int_{\mathbb{R}^{N}}\nabla
u_{n}\cdot\nabla(u_{n}-u)dx+\int_{\mathbb{R}^{N}}V(x)|u_{n}-u|^{2}dx\nonumber\\
& & \mbox{}-M(\|\nabla
u\|_{2}^{2})\int_{\mathbb{R}^{N}}\nabla
u\cdot\nabla(u_{n}-u)dx-\int_{\mathbb{R}^{N}}(f(x,u_{n})-f(x,u))(u_{n}-u)dx\nonumber\\
&=& \|\Delta (u_{n}-u)\|_{2}^{2}+M(\|\nabla
u_{n}\|_{2}^{2})\int_{\mathbb{R}^{N}}|\nabla(u_{n}-u)|^{2}dx+\int_{\mathbb{R}^{N}}V(x)|u_{n}-u|^{2}dx\nonumber\\
& & +\left(M(\|\nabla
u_{n}\|_{2}^{2})-M(\|\nabla
u\|_{2}^{2})\right)\int_{\mathbb{R}^{N}}\nabla
u\cdot\nabla(u_{n}-u)dx-\int_{\mathbb{R}^{N}}(f(x,u_{n})-f(x,u))(u_{n}-u)dx\nonumber\\
&\geq&
\|\Delta (u_{n}-u)\|_{2}^{2}+m_{0}\int_{\mathbb{R}^{N}}|\nabla(u_{n}-u)|^{2}dx+\int_{\mathbb{R}^{N}}V(x)|u_{n}-u|^{2}dx\nonumber\\
& & +\left(M(\|\nabla
u_{n}\|_{2}^{2})-M(\|\nabla
u\|_{2}^{2})\right)\int_{\mathbb{R}^{N}}\nabla
u\cdot\nabla(u_{n}-u)dx-\int_{\mathbb{R}^{N}}(f(x,u_{n})-f(x,u))(u_{n}-u)dx\nonumber\\
&\geq& \mbox{}
\min\{1,m_{0}\}\|u_{n}-u\|_{E}^{2}+
\left(M(\|\nabla
u_{n}\|_{2}^{2})-M(\|\nabla
u\|_{2}^{2})\right)\int_{\mathbb{R}^{N}}\nabla
u\cdot\nabla(u_{n}-u)dx\nonumber\\
&&-\int_{\mathbb{R}^{N}}(f(x,u_{n})-f(x,u))(u_{n}-u)dx.
\end{eqnarray*}
Define a linear functional $\mathcal{B}:E\rightarrow \mathbb{R}$ as
\begin{eqnarray*}
\mathcal{B}(\omega)=\int_{\mathbb{R}^{N}}\nabla u\cdot \nabla \omega dx.
\end{eqnarray*}
It can be deduced that $\mathcal{B}$ is continuous on $E$. Since
$u_{n}\rightharpoonup u$ in $E$, we obtain that
\begin{eqnarray*}
\int_{\mathbb{R}^{N}}\nabla
u\cdot\nabla(u_{n}-u)dx\rightarrow0\ \ \ \mbox{as}\ \ \
n\rightarrow\infty.
\end{eqnarray*}
Hence, by the boundedness of $\{u_{n}\}$ and the continuousness of
$M(t)$, we have
\begin{eqnarray*}
\left(M(\|\nabla
u_{n}\|_{2}^{2})-M(\|\nabla
u\|_{2}^{2})\right)\int_{\mathbb{R}^{N}}\nabla
u\cdot\nabla(u_{n}-u)dx\rightarrow0\ \ \ \mbox{as}\ \ \
n\rightarrow\infty,
\end{eqnarray*}
which implies that $\|u_{n}-u\|_{E}\rightarrow0$ as
$n\rightarrow\infty$. Hence $I$ satisfies the $(PS)$ condition.

From Lemmas 3.1-3.3 and the Mountain Pass Theorem, for any
$\lambda\in[0,\lambda_{1})$, we can obtain a critical point $u^{*}$ of
$I$ satisfying $I(u^{*})\geq\alpha$ and $I'(u^{*})=0$. The following
lemma tell us that there exists another nontrivial critical point of
$I$ corresponding to negative critical value.

\begin{lemma.} Suppose that (\ref{5}), $(V)$, $(M_{1})$, $(M_{2})$, $(f_{1})$, $(f_{2})$ and $(f_{6})$ hold, then
there exists a critical point of $I$ corresponding to negative
critical value for any $\lambda\in(0,\lambda_{1})$.
\end{lemma.}

{\bf Proof.} By Lemma 3.1, for any $\lambda\in(0,\lambda_{1})$, we can
see that there exits a local minimizer of $I$ in $B_{\varrho}$. The
following proof shows this minimizer is not zero. By $(f_{1})$, there
exists $\varsigma_{1}>0$ such that
\begin{eqnarray}
K(\bar{x},t)>\frac{1}{2}b_{0}|t|^{r_{0}}\label{34}
\end{eqnarray}
for all $x\in \Upsilon_{\varsigma_{1}}(\bar{x})$ and $t\in
\mathbb{R}$. Choosing $\varphi_{1}\in
C_{0}^{\infty}(\Upsilon_{\varsigma_{1}}(\bar{x}),\mathbb{R})\setminus\{0\}$,
by $(M_{2})$, $(f_{2})$, (\ref{34}) and $(f_{6})$, there exists
$M_{4}>0$ such that
\begin{eqnarray*}
I(\theta\varphi_{1})
&=&\frac{\theta^{2}}{2}\|\Delta \varphi_{1}\|_{2}^{2}+\frac{1}{2}\widehat{M}(\theta^{2}\|\nabla \varphi_{1}\|_{2}^{2})+\frac{\theta^{2}}{2}\int_{\Upsilon_{\varsigma_{1}}(\bar{x})}V(x)\varphi_{1}^{2}dx\\
&&-\lambda\int_{\Upsilon_{\varsigma_{1}}(\bar{x})}K(x,\theta \varphi_{1})dx-\int_{\Upsilon_{\varsigma_{1}}(\bar{x})}H(x,\theta \varphi_{1})dx\\
&\leq&\frac{\theta^{2}}{2}\|\Delta \varphi_{1}\|_{2}^{2}+\frac{\sigma_{2}}{2}\left(\theta^{4}\|\nabla \varphi_{1}\|_{2}^{4}+\theta^{2}\|\nabla \varphi_{1}\|_{2}^{2}\right)+\frac{\theta^{2}}{2}\int_{\Upsilon_{\varsigma_{1}}(\bar{x})}V(x)\varphi_{1}^{2}dx\\
&&-\frac{\lambda b_{0}\theta^{r_{0}}}{2}\int_{\Upsilon_{\varsigma_{1}}(\bar{x})}|\varphi_{1}|^{r_{0}}dx-M_{4}\int_{\Upsilon_{\varsigma_{1}}(\bar{x})}(\theta^{2}| \varphi_{1}|^{2}+\theta^{\zeta}| \varphi_{1}|^{\zeta} )dx\\
&<&0
\end{eqnarray*}
for $\theta>0$ small enough. Hence
\begin{eqnarray*}
-\infty<\inf\{I(u):u\in B_{\varrho}\}<0.
\end{eqnarray*}
Similar to the proof of Theorem 3.5 in \cite{20}, there exists
$u^{**}\in B_{\varrho}\setminus \partial B_{\varrho}$ such that
\begin{eqnarray*}
I(u^{**})=\inf_{u\in B_{\varrho}}I(u^{**})<0<\alpha\ \ \
\mbox{and}\ \ \ I^{\prime}(u^{**})=0.
\end{eqnarray*}
The proof of this lemma is finished.

\vspace{0.2cm} \textbf{Proof of Theorem \ref{th1.3}.} From Lemmas
3.1-3.4, we can see that problem (\ref{1}) possesses at least two
solutions for any $\lambda\in(0,\lambda_{1})$. \hfill$\Box$

\section{Proof of Theorem \ref{th1.4}}

In this section, we use Lemma \ref{leZhang} to obtain infinitely many
critical points of $I$. The following lemmas will show that $I$
satisfies the conditions of Lemma \ref{leZhang}.

\begin{lemma.}\label{le4.1}Suppose (\ref{5}), $(V)$, $(M_{1})$, $(M_{2})$, $(f_{2})$ and $(f_{3})$ hold, then I satisfies
$(C_{1})$. \end{lemma.}

\vspace{0.3cm}{\bf Proof.} Let $\{v_{j}\}_{j=1}^{\infty}$ be a
completely orthogonal basis of $E$ and
$X_{k}=\bigoplus_{j=1}^{k}S_{j}$, where $S_{j}=span\{v_{j}\}$. For any
$q\in[2,2^{*})$, we set
\begin{eqnarray}
\mathcal{P}_{k}(q)=\sup_{u\in X_{k}^{\bot},\|u\|_{E}=1}\|u\|_{q}.\label{22}
\end{eqnarray}
It follows from Lemma 2.10 in \cite{24} that
$\mathcal{P}_{k}(q)\rightarrow0$ as $k\rightarrow\infty$ for any
$q\in[2,2^{*})$. By $(f_{3})$, there exists $\rho_{1}\in(0,1)$ such
that
\begin{eqnarray}
|H(x,t)|\leq  \frac{t^{2}}{4C_{2}^{2}}\label{9}
\end{eqnarray}
for all $|t|\leq\rho_{1}$ and $x\in\mathbb{R}^{N}$ . Set
\begin{eqnarray}
\mathcal{H}_{k}=\frac{\lambda}{r_{1}}\mathcal{P}_{k}^{r_{1}}(r_{1}\beta_{1}^{*})\|b_{1}\|_{\beta_{1}}+\frac{\lambda}{r_{2}}\mathcal{P}_{k}^{r_{2}}(r_{2}\beta_{2}^{*})\|b_{2}\|_{\beta_{2}}+\frac{1}{4C_{2}^{2}}\mathcal{P}_{k}^{2}(2).\label{20}
\end{eqnarray}
Then there exists $k_{0}>0$ such that $
\mathcal{H}_{k}\leq\frac{1}{2}\min\left\{\frac{1}{2},\frac{m_{0}}{2}+\sigma_{1}\right\}\rho_{1}$
for all $k\geq k_{0}$.  Then for any $u\in
X_{k_{0}}^{\perp}\bigcap\partial B_{\rho_{1}}$, it follows from
$(M_{1})$, $(M_{2})$, $(f_{2})$, (\ref{22}) and (\ref{9}) that
\begin{eqnarray*}
I(u)&=&
\frac{1}{2}\|\Delta u\|_{2}^{2}+\frac{1}{2}\widehat{M}(\|\nabla u\|_{2}^{2})+\frac{1}{2}\int_{\mathbb{R}^{N}}V(x)u^{2}dx-\lambda\int_{\mathbb{R}^{N}}K(x,u)dx-\int_{\mathbb{R}^{N}}H(x,u)dx\nonumber\\
&\geq& \mbox{}
\min\left\{\frac{1}{2},\frac{m_{0}}{2}+\sigma_{1}\right\}\|u\|_{E}^{2}-\frac{\lambda}{r_{1}}\int_{\mathbb{R}^{N}}b_{1}(x)|u|^{r_{1}}dx-\frac{\lambda}{r_{2}}\int_{\mathbb{R}^{N}}b_{2}(x)|u|^{r_{2}}dx-\frac{1}{4C_{2}^{2}}\int_{\mathbb{R}^{N}}|u|^{2}dx\nonumber\\
&\geq& \mbox{}
\min\left\{\frac{1}{2},\frac{m_{0}}{2}+\sigma_{1}\right\}\rho_{1}^{2}-\frac{\lambda}{r_{1}}\mathcal{P}_{k}^{r_{1}}(r_{1}\beta_{1}^{*})\|b_{1}\|_{\beta_{1}}\rho_{1}^{r_{1}}-\frac{\lambda}{r_{2}}\mathcal{P}_{k}^{r_{2}}(r_{2}\beta_{2}^{*})\|b_{2}\|_{\beta_{2}}\rho_{1}^{r_{2}}-\frac{1}{4C_{2}^{2}}\mathcal{P}_{k}^{2}(2)\rho_{1}^{2}\nonumber\\
&\geq& \mbox{}
\min\left\{\frac{1}{2},\frac{m_{0}}{2}+\sigma_{1}\right\}\rho_{1}^{2}-\left(\frac{\lambda}{r_{1}}\mathcal{P}_{k}^{r_{1}}(r_{1}\beta_{1}^{*})\|b_{1}\|_{\beta_{1}}+\frac{\lambda}{r_{2}}\mathcal{P}_{k}^{r_{2}}(r_{2}\beta_{2}^{*})\|b_{2}\|_{\beta_{2}}+\frac{1}{4C_{2}^{2}}\mathcal{P}_{k}^{2}(2)\right)\rho_{1}\\
&\geq&\frac{1}{2}\min\left\{\frac{1}{2},\frac{m_{0}}{2}+\sigma_{1}\right\}\rho_{1}^{2}.
\end{eqnarray*}
We finish the proof of this lemma.

\begin{lemma.}\label{le4.2} Suppose (\ref{5}), $(V)$, $(M_{2})$, $(f_{2})$ and $(f_{4})$ hold, then $I$ satisfies $(C_{2})$. \end{lemma.}

{\bf Proof.} Set $\tilde{X}_{m}=\bigoplus_{j=1}^{m}S_{j}$, where
$S_{j}$ is defined in Lemma 4.1. For any $u\in
\tilde{X}_{m}\setminus\{0\}$ and $\vartheta>0$, set
\begin{eqnarray*}
\Gamma_{\vartheta}(u)=\{x\in \mathbb{R}:\ |u|\geq\vartheta \|u\|_{E}\}.
\end{eqnarray*}
Similar to \cite{25}, there exists $\vartheta_{0}>0$ such that
\begin{eqnarray*}
meas\left(\Gamma_{\vartheta_{0}}(u)\right)\geq\vartheta_{0}\label{21}
\end{eqnarray*}
for all $u\in \tilde{X}_{m}\setminus\{0\}$. Then there exists
$\kappa>0$ such that
\begin{eqnarray}
meas\left(\Lambda_{\vartheta_{0}}(u)\right)\geq\frac{3}{4}\vartheta_{0},\label{18}
\end{eqnarray}
for all $u\in \tilde{X}_{m}\setminus\{0\}$, where
$\Lambda_{\vartheta_{0}}(u)=\Gamma_{\vartheta_{0}}(u)\bigcap
\Upsilon_{\kappa}(0)$. It follows from $(f_{4})$ that there exists
$L_{2}>0$ such that
\begin{eqnarray}
H(x,u)\geq \frac{\sigma_{2}}{\vartheta_{0}^{5}}|u|^{4}\geq
\frac{\sigma_{2}}{\vartheta_{0}}\|u\|_{E}^{4}\label{21}
\end{eqnarray}
for all $u\in \tilde{X}_{m}$ and $x\in \Lambda_{\vartheta_{0}}(u)$
with $\|u\|_{E}\geq L_{2}$. We can choose $\varsigma_{m}>L_{2}$, then
for any $u\in\tilde{X}_{m}\setminus B_{\varsigma_m}$, it follows from
(\ref{3}), $(M_{2})$, (\ref{27}), (\ref{6}), (\ref{18}) and (\ref{21})
that
\begin{eqnarray*}
I(u)&=&\frac{1}{2}\|\Delta u\|_{2}^{2}+\frac{1}{2}\widehat{M}(\|\nabla u\|_{2}^{2})+\frac{1}{2}\int_{\mathbb{R}^{N}}V(x)u^{2}dx-\lambda\int_{\mathbb{R}^{N}}K(x,u)dx-\int_{\mathbb{R}^{N}}H(x,u)dx\nonumber\\
&\leq& \frac{1}{2}\|\Delta u\|_{2}^{2}+\frac{\sigma_{2}}{2}(\|\nabla u\|_{2}^{4}+\|\nabla u\|_{2}^{2})+\frac{1}{2}\int_{\mathbb{R}^{N}}V(x)u^{2}dx-\lambda\int_{\mathbb{R}^{N}}K(x,u)dx-\int_{\mathbb{R}^{N}}H(x,u)dx{\nonumber}\\
&\leq& \frac{1+\sigma_{2}}{2}\|u\|_{E}^{2}+\frac{\sigma_{2}}{2}\|u\|_{E}^{4}+\lambda M_{1}(\|u\|_{E}^{r_{1}}+\|u\|_{E}^{r_{2}})-\int_{\Lambda_{\vartheta_{0}}(u)}H(x,u)dx+M_{3}\int_{\mathbb{R}^{N}\setminus\Lambda_{\vartheta_{0}}(u)}|u|^{2}dx{\nonumber}\\
&\leq& \frac{1+\sigma_{2}}{2}\|u\|_{E}^{2}+\frac{\sigma_{2}}{2}\|u\|_{E}^{4}+\lambda M_{1}(\|u\|_{E}^{r_{1}}+\|u\|_{E}^{r_{2}})-\frac{\sigma_{2}}{\vartheta_{0}}meas(\Lambda_{\vartheta_{0}}(u))\|u\|_{E}^{4}+M_{3}C_{2}^{2}\|u\|_{E}^{2}{\nonumber}\\
&\leq& -\frac{\sigma_{2}}{4}\|u\|_{E}^{4}+\left(\frac{1+\sigma_{2}}{2}+M_{3}C_{2}^{2}\right)\|u\|_{E}^{2}+\lambda M_{1}(\|u\|_{E}^{r_{1}}+\|u\|_{E}^{r_{2}}).
\end{eqnarray*}
Then there exists $r_m>\xi$ such that $I(u_{m})\leq0$ for all
$u\in\tilde{X}_{m}\setminus B_{r_m}$, which proves this lemma.

\begin{lemma.}\label{le4.3} Suppose the conditions of
Theorem \ref{th1.4} hold, then $I$ satisfies the $(PS)^{*}$ condition.
\end{lemma.}

{\bf Proof.} The proof is similar to Lemma \ref{le3.3}, we omit it
here.

\vspace{0.2cm} \textbf{Proof of Theorem \ref{th1.4}.} By Lemmas
\ref{le4.1}-\ref{le4.3} and Lemma \ref{leZhang}, $I$ possesses
infinitely many distinct critical points corresponding to positive
critical values. \hfill$\Box$

\section{Acknowledgments}

The authors are grateful to the referees for the helpful comments
which improve the writing of the paper. This paper was finished when
D.-L. Wu was visiting Utah State University with the support of China
Scholarship Council(No.201708515186); he is grateful to the members in
the Department of Mathematics and Statistics at Utah State University
for their invitation and hospitality.

\def\refname{References}

\end{document}